\documentclass[12pt]{article}
\usepackage{amscd,amsfonts,amssymb,amsmath,latexsym,array,hhline}

\renewcommand{\bar}{\overline}
\renewcommand{\tilde}{\widetilde}
\renewcommand{\hat}{\widehat}

\newcommand{\LL}{{\mathbb L}}

\newcommand{\NN}{{\mathbb N}}
\newcommand{\RR}{{\mathbb R}}

 \numberwithin{equation}{section}

 \renewcommand{\section}[1]{\refstepcounter{section}
                            {\noindent\large\bf\thesection. #1}}

\makeatletter

\renewcommand{\section}{\@startsection{section}{1}{0pt}{30pt}{6pt}{\large\bf}}
\def\dot{\hspace{-16pt}.\hspace{-2pt} }

\renewcommand{\@makefnmark}{}
\renewcommand{\@cite}[2]{[{#1\if@tempswa ; #2\fi}]}

\makeatother

\begin{document}

\title{\bf \LARGE Bayesian quickest detection problems \\
                  for some diffusion processes}

\author{Pavel V. Gapeev$^{*}$ \and Albert N. Shiryaev$^{**}$}
\date{}
\maketitle


\begin{abstract}
        We study the Bayesian problems of detecting a change in the drift
        rate of an observable diffusion process with linear and exponential
        penalty costs for a detection delay.
        The optimal times of alarms are found as the first times at which the weighted
        likelihood ratios hit stochastic boundaries depending on the current observations.
        The proof is based on the reduction of the initial problems into
        appropriate three-dimensional optimal stopping problems and the
        analysis of the associated parabolic-type free-boundary problems.
        We provide closed form estimates for the value functions
        and the boundaries, 
        under certain nontrivial relations between the coefficients
        of the observable diffusion.
\end{abstract}
%







\footnotetext{$^{*}$\, London School of Economics, Department of
Mathematics, Houghton Street, London WC2A 2AE, United Kingdom;
e-mail: p.v.gapeev{\char'100}lse.ac.uk} \footnotetext{$^{**}$\,
Steklov Institute of Mathematics, Russian Academy of Sciences,
Gubkina Street 8, Moscow 119991, Russia; e-mail:
albertsh{\char'100}mi.ras.ru}



\section{\dot Introduction}

     The problem of quickest disorder detection for an observable
     diffusion process seeks to determine a stopping time of alarm
     $\tau$ which is as close as possible to the unknown time
     of {\it disorder} (or change-point) $\theta$ at which the local drift rate 
     of the process changes from $\mu_0(\cdot)$ to $\mu_1(\cdot)$.
     In the classical Bayesian formulation, it is assumed that the random time $\theta$
     takes the value $0$ with probability $\pi$ and is exponentially distributed with parameter
     $\lambda > 0$ given that $\theta > 0$.
     An optimality criterion was proposed in \cite{S61b}-\cite{S63} for the time of alarm
     to minimize a linear combination of the false alarm probability and the expected time
     delay in detecting the disorder correctly, for sequences of i.i.d. observations.
     An explicit solution of the problem of detecting a change in the constant drift rate
     of an observable Wiener process with the same optimality criterion was derived in
     \cite{S65}-\cite{SCyb}. The appropriate optimal stopping problem for the
     posterior probability of the occurrence of disorder was reduced to the
     associated free-boundary problem for an ordinary differential operator
     (see also \cite[Chapter~IV, Section~4]{S} or \cite[Chapter~VI, Section~22]{PSbook}).
     A finite time horizon version of the Wiener disorder problem
     was studied in \cite{GP2}.

     The idea of replacing the initial average time delay by a certain non-additive detection delay penalty criterion
     was originally introduced in \cite{S64}.
     The resulting Bayesian risk function was expressed through the current state of a multi-dimensional
     Markovian sufficient statistic having state space components which are different from the posterior
     probability. Such a process contained all the necessary information to determine the
     structure of the optimal time of alarm (see also more recent works \cite{S07}, \cite{SZ} and \cite{Dayanik}).
     In the case of exponential penalty costs for a delay, it was observed by Poor \cite{Poor} that
     the weighted likelihood ratio process turns out to be a one-dimensional Markovian sufficient
     statistic, for sequences of i.i.d. observations. This idea was taken further by Beibel
     \cite{Beibel}, who solved the corresponding problem of detecting a change in the drift
     rate of an observable Wiener process as a generalized parking problem. Bayraktar and
     Dayanik \cite{BD} recognized the same property
     from the structure of the
     ordinary differential-difference equation in the free-boundary problem associated with
     the Bayesian problem of detecting a change in the constant intensity rate of an observable
     Poisson process. Some other formulations of the problem for the case of detecting a change
     in the arrival rate of a Poisson process, leading to the appearance of essentially
     multi-dimensional Markovian sufficient statistics, were studied by Bayraktar,
     Dayanik, and Karatzas \cite{BDK1}-\cite{BDK2}. Extensive overviews of these and other related quickest
     sequential change-point detection methods were provided in the monographs \cite{S02} and \cite{PH}.

     In the present paper, we study the Bayesian quickest disorder detection problems
     for observable diffusions with linear and exponential delay penalty costs.
     We reduce the initial problems to extended optimal stopping problems for
     three-dimensional Markov diffusion processes, having the posterior
     probability, weighted likelihood ratio, and the observations
     as their state space components. We show that the optimal stopping times
     are expressed as the first times at which the weighted likelihood ratio
     processes hit stochastic boundaries depending on the current state of
     the observation process only. We verify that the value
     functions and the optimal stopping boundaries are characterized
     by means of the associated free-boundary problem for a second-order
     partial differential operator. The latter turns out to be of parabolic
     type, because the observation process is a one-dimensional diffusion.
     We also derive closed form estimates for the value functions and the
     boundaries for a special nontrivial subclass of observable diffusions.
     The Bayesian sequential testing problem for such
     processes was recently solved in \cite{GS1}.
     Another related problem of transient signal detection and
     identification of two-sided changes in the drift rates of
     observable diffusion processes was considered by Pospisil,
     Vecer and Hadjiliadis \cite{PVH}.

     The paper is organized as follows.
     In Section 2, we formulate the Bayesian quickest disorder detection problem for observable diffusion
     processes with linear and exponential delay penalty costs and construct the associated
     multi-dimensional optimal stopping problem.
     In Sections 3 and 4, we present the associated free-boundary problem and
     reduce the resulting parabolic-type partial differential
     operator to the normal form, which is amenable for further considerations.
     Applying the change-of-variable formula with local
     time on surfaces, obtained by Peskir \cite{Pe1a}, we verify that the solution
     of the free-boundary problem, which satisfies certain additional conditions,
     provides the solution of the initial optimal stopping problem.
     We derive closed form estimates for the value function
     and the boundary, which are uniquely determined as solutions of
     ordinary differential equations,
     under certain nontrivial relations between the coefficients
     of the observable diffusion. The main results are stated in Theorems 3.4 and 4.2.

  \section{\dot Preliminaries}


  In this section, we give the Bayesian formulation of the problem
  (see \cite[Chapter~IV, Section~4]{S} or \cite[Chapter~VI, Section~22]{PSbook} for the case of Wiener
  processes) in which it is assumed that one observes a sample path of the diffusion process
  $X=(X_t)_{t \ge 0}$ with the drift rate changing from $\mu_0(\cdot)$
  to $\mu_1(\cdot)$ at some random time $\theta$ taking the value $0$
  with probability $\pi$ and being exponentially distributed with
  parameter $\lambda > 0$ under $\theta > 0$.

  \vspace{6pt}

  2.1. (Formulation of the problem.) 
  Suppose that, on a probability space $(\Omega, {\cal F}, {P}_{\pi})$, there exists
  a standard Brownian motion $B=(B_t)_{t \ge 0}$ independent of a nonnegative random
  variable $\theta$ such that ${P}_{\pi}(\theta=0) = \pi$ and
  ${P}_{\pi}(\theta > t \, | \, \theta>0) = e^{-\lambda t}$,
  for all $t \ge 0$ and some $\lambda > 0$ fixed.
  Let $X=(X_t)_{t \ge 0}$ be a continuous process solving
  the stochastic differential equation:
  \begin{equation}
  \label{X3}
  dX_t = \big( \mu_0(X_t) + I(\theta \le t) (\mu_1(X_t)-\mu_0(X_t)) \big) \, dt +
  \sigma(X_t) \, dB_t
  \end{equation}
  with $X_0=x$, where $\mu_i(x)$, $i = 0, 1$, and $\sigma(x) > 0$ are some continuously
  differentiable functions on $(0, \infty)$, satisfying the 
  conditions:
  \begin{equation}
  \label{musig5}
  |\mu_i(x)| + |\sigma(x)| \le K \, (1 + |x|) \quad \text{and} \quad
  0 < \bigg| \frac{\mu_1(x)-\mu_0(x)}{\sigma(x)} \bigg| \le K
  \end{equation}
  for all $x > 0$ and some $K > 0$ fixed.
  In order to facilitate the considerations of the examples below,
  we assume the state space of the process $X$ to be the positive
  half line $(0, \infty)$.
  It thus follows from \cite[Theorem~4.6]{LS} that the equation
  in (\ref{X3}) admits a unique strong solution under $\theta=s$, and hence,
  ${P}_{\pi}(X \in \: \cdot \;| \, \theta=s \,)=P^s(X \in \: \cdot \;)$
  is the distribution law of a time-homogeneous diffusion process
  started at some $x > 0$,
  with diffusion coefficient $\sigma(x)$ and the drift rate
  changing from $\mu_0(x)$ to $\mu_1(x)$ at time $s \in [0, \infty]$.
  In this case, we may conclude that the probability measure $P_{\pi}$
  has the structure:
  \begin{equation}
  \label{ppi3}
  {P}_{\pi}(X \in \: \cdot \;) = \pi P^0(X \in \: \cdot \;) + (1-\pi) \int_0^{\infty}
  P^s(X \in \: \cdot \;) \, \lambda e^{-\lambda s} \, ds
  \end{equation}
  for any $\pi \in [0, 1)$ fixed.

  Based upon the continuous observation of the process $X$, our task is to find among
  the stopping times $\tau$ of $X$ (i.e. stopping times with respect to the natural filtration
  ${\cal F}_t=\sigma (X_s \, | \, 0 \le s \le t)$ of the process $X$)
  an optimal time at which an {\it alarm} should be sounded {\it as close as possible}
  to the unobservable time of {\it disorder} $\theta$.
  More precisely, the {\it Bayesian quickest detection problem}
  consists of computing the Bayesian risk function:
  \begin{equation}
  \label{V2}
  {V}(\pi) = \inf_{\tau}
  \Big( P_{\pi}(\tau<\theta) + E_{\pi} [F(\tau-\theta) I(\tau \ge \theta)] \Big)
  \end{equation}
  and finding the optimal stopping time, called the $\pi$-Bayesian time,
  at which the infimum is attained in $(\ref{V2})$.
  Here $P_{\pi}(\tau < \theta)$ is the probability of a {\it false alarm},
  and $E_{\pi}[F(\tau-\theta)I(\tau \ge \theta)]$ is the expected {\it costs of delay}
  in detecting of the disorder correctly (i.e. when $\tau \ge \theta$), where
  the {\it delay penalty function} $F(t)$ satisfies the conditions $F(t) \ge 0$ for $t \ge 0$,
  and $F(t) = 0$ for $t \le 0$. We will further assume that either $F(t) = c t$ or
  $F(t) = c(e^{\alpha t} - 1)$ holds in (\ref{V2}) for all $t \ge 0$.

   \vspace{6pt}

   {\bf Remark 2.1.} It was shown in
   \cite{S64}, 
   \cite{S07} and \cite{Dayanik} that,
   when the Laplace transforms of delay penalty functions are of rational structure,
   there exist finite-dimensional
   processes called Markovian {\it sufficient statistics} in the corresponding Bayesian
   quickest detection problems. Such (time-homogeneous strong) Markov processes containing all the necessary
   information to determine the optimal stopping times (see \cite[Chapter~II, Section~15]{S}
   for an extensive discussion of this notion).
   For example, the function $F(t) = c t^{\delta}$ for $t \ge 0$,
   with some $c, \delta > 0$, $\delta \in \NN$, is of such type,
   while the assumption $\delta \notin \NN$ leads to the appearance of
   an infinite-dimensional Markovian sufficient statistic in that case.

  \vspace{6pt}

  2.2. (Likelihood ratio and posterior probability.) 
In order to derive Markovian sufficient statistics for the problem of (\ref{V2}), for the cases
of linear and exponential delay penalty functions indicated above,
let us define the {\it posterior probability} process $(\pi_t)_{t \ge 0}$
by $\pi_t = P(\theta \le t \, | \, {\cal F}_t)$ for $t \ge 0$.
  Taking into account the fact that the probability measure $P^s$ is equivalent
  to $P_{\pi}$ on ${\cal F}_t$ by construction, for any $s \in [0, \infty]$,
  using Bayes' formula (see, e.g. \cite[Theorem~7.23]{LS}), we get that
  $(\pi_t)_{t \ge 0}$ admits the representation:
  \begin{equation}
  \label{V2b}
  \pi_t = \pi \, \frac{d(P^0 \, | \, {\cal
  F}_t)}{d(P_{\pi} \, | \, {\cal F}_t)} + (1-\pi) \int_0^t \frac{d(P^s \, | \, {\cal
  F}_t)}{d(P_{\pi} \, | \, {\cal F}_t)} \, \lambda e^{-\lambda s} \, ds.
  \end{equation}
  Moreover, since the measure $P^u$ coincides with $P^t$ on ${\cal F}_t$,
  for all $0 \le t \le u$, we see that:
  \begin{equation}
  \label{V2c}
  1 - \pi_t = (1-\pi) \int_t^{\infty}
  \frac{d(P^u \, | \, {\cal F}_t)}{d(P_{\pi} \, | \, {\cal F}_t)} \, \lambda e^{-\lambda u} \, du
  = (1-\pi) \, e^{-\lambda t} \, \frac{d(P^t \, | \, {\cal F}_t)}{d(P_{\pi} \, | \, {\cal F}_t)}
  \end{equation}
  is satisfied. By means of Girsanov's theorem for diffusion processes
  (see, e.g. \cite[Theorem~7.19]{LS}), it follows from the structure of the observation
  process $X$ in (\ref{X3}) that the {\it likelihood ratio} process $L = (L_t)_{t \ge 0}$ defined by:
  \begin{equation}
  \label{V2e}
  L_t = \frac{d(P^0 \, | \, {\cal F}_t)}{d(P^t \, | \, {\cal F}_t)}
  \equiv \frac{d(P^0 \, | \, {\cal F}_t)}{d(P^{\infty} \, | \, {\cal F}_t)}
  \end{equation}
  admits the representation:
  \begin{equation}
  \label{V2f}
  L_t = \exp \left( \int_0^t \frac{\mu_1(X_s)-\mu_0(X_s)}{\sigma^2(X_s)} \,
  dX_s - \frac{1}{2} \int_0^t \frac{\mu_1^2(X_s)-\mu_0^2(X_s)}{\sigma^2(X_s)} \,
  ds \right).
  \end{equation}
  Hence, the expressions in (\ref{V2e}) and (\ref{V2f}) yield that the properties:
  \begin{equation}
  \label{V2cc}
  \frac{d(P^s \, | \, {\cal F}_t)}{d(P_{\pi} \, | \, {\cal F}_t)} \,
  \frac{d(P_{\pi} \, | \, {\cal F}_t)}{d(P^t \, | \, {\cal F}_t)}
  = \frac{d(P^s \, | \, {\cal F}_t)}{d(P^0 \, | \, {\cal F}_t)} \,
  \frac{d(P^0 \, | \, {\cal F}_t)}{d(P^t \, | \, {\cal F}_t)}
  = \frac{d(P^s \, | \, {\cal F}_s)}{d(P^0 \, | \, {\cal F}_s)} \,
  \frac{d(P^0 \, | \, {\cal F}_t)}{d(P^t \, | \, {\cal F}_t)} \equiv \frac{L_t}{L_s}
  \end{equation}
  hold for each $0 \le s \le t$. We therefore obtain from the representations
  in (\ref{V2b}) and (\ref{V2c}) that the {\it weighted likelihood ratio}
  process $(\varphi_t)_{t \ge 0}$ defined by ${\varphi}_t = {\pi_t}/({1-\pi_t})$
  has the form:
  \begin{equation}
  \label{f2g}
  {\varphi}_t = e^{\lambda t} L_t \left(\frac{\pi}{1-\pi} +
  \int_0^t \frac{\lambda e^{- \lambda s}}{L_s} \, ds \right).
  \end{equation}


\vspace{3pt}

  2.3. (Stochastic differential equations.)
  Applying It\^o's formula
  (see, e.g. \cite[Chapter~IV, Theorem~4.4]{LS} or \cite[Chapter~IV, Theorem~3.3]{RY})
  to the expression in (\ref{V2f}), we get that the process $L$ admits the representation:
  \begin{equation}
  \label{L22}
  dL_t = 
  \frac{\mu_1(X_t)-\mu_0(X_t)}{\sigma^2(X_t)} \, L_{t} \, (dX_t - \mu_0(X_t) \, dt)
  \end{equation}
  with $L_0 = 1$. Then, using the integration by parts formula, we see that
  the process $({\varphi}_t)_{t \ge 0}$ from (\ref{f2g}) solves the stochastic
  differential equation:
   \begin{equation}
   \label{phi3d}
   d {\varphi}_t = \bigg( \lambda (1 + {\varphi}_t) + \bigg( \frac{\mu_1(X_t)-\mu_0(X_t)}{\sigma(X_t)} \bigg)^2 \,
   \frac{{\varphi}^2_t}{1+{\varphi}_t} \bigg) \, dt + \frac{\mu_1(X_t)-\mu_0(X_t)}{\sigma(X_t)} \, {\varphi}_t \, d{\bar B}_t
   \end{equation}
   with $\varphi_0 = \varphi \equiv \pi/(1-\pi)$. Hence, using It{\^o}'s formula again, we obtain that
   the process $(\pi_t)_{t \ge 0}$ admits the representation:
   \begin{equation}
   \label{pi3d}
   d\pi_t = \lambda (1-\pi_t) \, dt +
   \frac{\mu_1(X_t)-\mu_0(X_t)}{\sigma(X_t)} \, \pi_t(1-\pi_t) \, d{\bar B}_t 
   \end{equation}
   with $\pi_0 = \pi$. Here, the {\it innovation} process ${\bar B}=({\bar B}_t)_{t \ge 0}$ defined by:
       \begin{equation}
       \label{W2}
       {\bar B}_t = \int_0^t \frac{dX_s}{\sigma(X_s)} -
       \int_0^t \left( \frac{\mu_0(X_s)}{\sigma(X_s)} +
       \pi_s \, \frac{\mu_1(X_s)-\mu_0(X_s)}{\sigma(X_s)} \right) ds
       \end{equation}
       is a standard Brownian motion under the measure $P_{\pi}$,
       with respect to the filtration $({\cal F}_t)_{t \ge 0}$,
       according to P. L{\'e}vy's characterization theorem
       (see, e.g. \cite[Theorem~4.1]{LS} or \cite[Chapter~IV, Theorem~3.6]{RY}).
       It thus follows from (\ref{W2}) that the process $X$ admits the representation:
       \begin{equation}
       \label{X3a}
       dX_t = \big( \mu_0(X_t) + \pi_t \, (\mu_1(X_t)-\mu_0(X_t)) \big) \, dt
       + \sigma(X_t) \, d{\bar B}_t 
       \end{equation}
       with $X_0=x$. Taking into account the assumptions in (\ref{musig5}),
       we may conclude by virtue of Remark to \cite[Chapter~IV, Theorem~4.6]{LS}
       (see also \cite[Chapter~V, Theorem~5.2.1]{O}) that the processes
       $(\pi_t, X_t)_{t \ge 0}$ and $(\varphi_t, X_t)_{t \ge 0}$ turn out
       to be unique strong solutions of the corresponding systems of
       stochastic differential equations in (\ref{phi3d}), (\ref{pi3d}), and (\ref{X3a}).
       According to \cite[Chapter~VII, Theorem~7.2.4]{O},
       such processes have the (time-homogeneous strong) Markov property
       with respect to its natural filtration, which inherently coincides
       with $({\cal F}_t)_{t \ge 0}$.

\vspace{6pt}

2.4. (Some examples.) Let us now present some expressions for the Bayesian risk functions
and the appropriate Markovian sufficient statistics in the corresponding quickest
disorder detection problems for observable diffusion processes.



\vspace{6pt}

{\bf Example 2.2.} Assume that we have $F(t) = c t$ with 
some $c > 0$ fixed (see \cite{S65}, \cite{SCyb},
\cite[Chapter~IV]{S}, and \cite[Chapter~VI, Section~22]{PSbook}). It is then shown by means of standard arguments from
\cite[Chapter~IV, Section~3]{S} that the Bayesian risk function ${V}(\pi)$ in (\ref{V2}) admits the representation:
\begin{equation}
\label{V30}
V'(\pi, \varphi, x) = \inf_{\tau} E_{\pi, \varphi, x} \left[1 - \pi_{\tau} + \int_0^{\tau} (1 - \pi_t) \, c \varphi_t \, dt \right]
\end{equation}
where the infimum is taken over all stopping times $\tau$ such that $E_{\pi, \varphi, x} \tau < \infty$ holds.
Here, $P_{\pi, \varphi, x}$ is a measure of the diffusion process $(\pi_t, \varphi_t, X_t)_{t \ge 0}$, started at some
$(\pi, \varphi, x) \in [0, 1) \times [0, \infty) \times (0, \infty)$ and solving the equations in (\ref{phi3d}),
(\ref{pi3d}), and (\ref{X3a}), which is a Markovian sufficient statistic in the problem.

\vspace{6pt}

  {\bf Example 2.3.} Assume now that $F(t) = c (e^{\alpha t} - 1)$ with
  some $c, \alpha > 0$
  fixed (see \cite[Example~4]{S64}, \cite{Poor}, \cite{Beibel}, and \cite{BD}).
  It can be shown following the schema of arguments from \cite{BD} that the
  Bayesian risk function ${V}(\pi)$ in (\ref{V2}) admits the representation: 
  \begin{equation}
  \label{V3b}
  {V}_*(\pi, \phi, x) = \inf_{\tau} E_{\pi, \phi, x} \left[1 - \pi_{\tau} + \int_0^{\tau} (1-\pi_t) \, c \alpha \phi_t \, dt \right]
  \end{equation}
  where the infimum is taken over all stopping times $\tau$ such that the integral above has a finite
  expectation, so that $E_{\pi, \phi, x} \tau < \infty$ holds.
   Here, the {\it weighted likelihood ratio} process $(\phi_t)_{t \ge 0}$ defined by:
   \begin{equation}
   \label{Phi2a}
   \phi_t = e^{(\alpha + \lambda) t} {L_t} \left(\frac{\pi}{1-\pi}+
   \int_0^t \frac{\lambda e^{- (\alpha + \lambda) s}}{L_s} \, ds \right)
   \end{equation}
   solves the stochastic differential equation:
   \begin{equation}
   \label{Phi2d}
   d \phi_t = \bigg( \lambda + (\lambda + \alpha) \, \phi_t + \bigg( \frac{\mu_1(X_t)-\mu_0(X_t)}{\sigma(X_t)} \bigg)^2 \, \pi_t \, \phi_t
   \bigg) \, dt + \frac{\mu_1(X_t)-\mu_0(X_t)}{\sigma(X_t)} \, \phi_t \, d{\bar B}_t
   \end{equation}
   with $\phi_0 = \phi \equiv \pi/(1-\pi)$.
   In this case, $P_{\pi, \phi, x}$ is a measure of the diffusion process $(\pi_t, \phi_t, X_t)_{t \ge 0}$,
   started at some $(\pi, \phi, x) \in [0, 1) \times [0, \infty) \times (0, \infty)$ and solving the equations
   in (\ref{pi3d}), (\ref{Phi2d}), and (\ref{X3a}), which is a Markovian sufficient statistic in the problem.

  \section{\dot The case of exponential delay penalty costs}


In this section, we formulate and prove the main assertions of the
paper, which are related to the quickest detection problem with exponential
delay penalty costs 
of Example 2.3 above.


        \vspace{6pt}


        3.1. By means of the results of general theory of optimal stopping
        (see, e.g. \cite[Chapter~III]{S} or \cite[Chapter~I, Section~2.1]{PSbook}), it follows
        from the structure of the reward functional in (\ref{V3b}) that the optimal stopping
        time is given by:
        \begin{equation}
        \label{tau0} \tau_* = \inf \{t \ge 0 \, | \, V_*(\pi_t, \phi_t, X_t) =  1 - \pi_t \}
        \end{equation}
        whenever the corresponding integral there is of finite expectation, so that $E_{\pi, \phi, x} \tau_* < \infty$ holds.
        In order to specify the structure of the stopping time in (\ref{tau0}),
        we follow the arguments from \cite[Subsection~2.5]{GP2} and use It{\^o}'s formula to get:
        \begin{equation}
        \label{GIto3}
        1 - \pi_{t} = 1 - \pi - \int_0^t \lambda \, (1 - \pi_{s}) \, ds + N_t
        \end{equation}
        where the process $N = (N_t)_{t \ge 0}$ defined by:
        \begin{equation}
        \label{N3}
        N_t = - \int_0^t
        \frac{\mu_1(X_s)-\mu_0(X_s)}{\sigma(X_s)} \, \pi_s (1-\pi_s) \, d{\bar B}_s
        \end{equation}
        is a continuous local martingale under $P_{\pi, \phi, x}$. It follows directly from (\ref{GIto3})
        that the process $(N_{\tau \wedge t})_{t \ge 0}$ is a uniformly integrable martingale
        for any stopping time $\tau$ satisfying $E_{\pi, \phi, x} \tau < \infty$.
        Then, applying Doob's optional sampling theorem (see, e.g.
        \cite[Theorem~3.6]{LS} or \cite[Chapter~II, Theorem~3.2]{RY}),
        we get from the expression in (\ref{GIto3}) that:
        \begin{equation}
        \label{2.13}
        E_{\pi, \phi, x} \left[ 1 - \pi_{\tau} + \int_0^{\tau}
        (1-\pi_{t}) \, c \alpha \phi_{t} \, dt \right]
        = 1 - \pi + E_{\pi, \phi, x} \int_0^{\tau} (1-\pi_{t})
        \, (c \alpha \phi_t - \lambda) \, dt
        \end{equation}
        holds for all $(\pi, \phi, x) \in [0, 1) \times [0, \infty) \times (0, \infty)$
        and any $\tau$ such that $E_{\pi, \phi, x} \tau < \infty$.
        Taking into account the structure of the reward in 
        (\ref{V3b}), it is seen from (\ref{2.13}) that it is never optimal to stop when
        $\phi_t < \lambda/(c \alpha)$ for any $t \ge 0$.
        This shows that all the points $(\pi, \phi, x)$
        such that $\phi < \lambda/(c \alpha)$ belong to the continuation region:
        \begin{equation}
        \label{C3aa}
        C_* = \{ (\pi, \phi, x) \in [0, 1) \times [0, \infty) \times (0, \infty) \, | \, V_*(\pi, \phi, x) <  1 - \pi \}.
        \end{equation}

%

        \vspace{3pt}

        3.2. In order to describe the structure of the set in (\ref{C3aa}),
        let us fix some $(\pi, \phi, x) \in C_*$ and denote by $\tau_*=\tau_*(\pi, \phi, x)$
        the optimal stopping time in the problem of (\ref{V3b}).
        Then, by means of the general optimal stopping theory
        for Markov processes (see, e.g. \cite[Chapter~III]{S}
        or \cite[Chapter~I, Section~2.2]{PSbook}), we conclude that:
        \begin{equation}
        \label{Vsigma3c}
        V_*(\pi, \phi, x) 
        = E_{\pi, \phi, x}
        \left[ 1 - \pi_{\tau_*} +
        \int_0^{\tau_*} (1-\pi_{t}) \, c \alpha \phi_{t} \, dt \right] 
        < 1 - \pi
        \end{equation}
        holds. Hence, taking any $\phi'$ such that $\phi' < \phi$
        and using the explicit expression for the process $(\phi_t)_{t \ge 0}$ through its
        starting point $\phi \equiv \pi/(1-\pi)$ in (\ref{Phi2a}), we obtain from (\ref{V3b})
        that the inequalities:
        \begin{align}
        \label{Vsigma3c'}
        V_*(\pi, \phi', x) 
        &\le E_{\pi, \phi', x} \left[ 1 - \pi_{\tau_*}
        + \int_0^{\tau_*} (1-\pi_{t}) \, c \alpha \phi_{t} \, dt \right] 
        \\
        \notag
        &\le E_{\pi, \phi, x} \left[ 1 - \pi_{\tau_*} + \int_0^{\tau_*} (1-\pi_{t}) \,
        c \alpha \phi_{t} \, dt \right] 
        \end{align}
        are satisfied. Thus, by virtue of the inequality in (\ref{Vsigma3c}), we see that $(\pi, \phi', x) \in C_*$.
        Taking into account the multiplicative structure of the integrand in (\ref{V3b}),
        we can therefore extend the approach used in \cite{Poor}, \cite{Beibel}, and \cite{BD},
        and further assume that there exists a function $g_*(x)$ such that $0 < \lambda/(c \alpha) \le g_*(x)$ for $x > 0$,
        and the continuation region in (\ref{C3aa}) for the optimal stopping problem of (\ref{V3b}) takes the form:
        \begin{equation}
        \label{C5}
        C_* = \{(\pi, \phi, x) \in [0, 1) \times [0, \infty) \times (0, \infty) \, | \, \phi < g_*(x) \}
        \end{equation}
        so that the corresponding stopping region is the closure of the set:
        \begin{equation}
        \label{D5}
        D_* = \{(\pi, \phi, x) \in [0, 1) \times [0, \infty) \times (0, \infty) \, | \, \phi > g_*(x) \}.
        \end{equation}

    \vspace{3pt}

        3.3. In order to characterize the behavior of the boundary $g_*(x)$ in (\ref{C5})-(\ref{D5}),
        we observe from the equation in (\ref{Phi2d}) that the expression:
        \begin{equation}
        \label{GIto4}
        \int_0^{\tau_*} ({1-\pi_t}) \, \phi_t \, dt = \int_0^{\tau_*} \frac{1-\pi_t}{\lambda + \alpha} \, d\phi_t
        - \int_0^{\tau_*} \frac{1-\pi_t}{\lambda + \alpha} \, (\lambda + \rho(X_t) \, \pi_t \, \phi_t) \, dt + N^*_{\tau_*}
        \end{equation}
        holds for the optimal stopping time $\tau_*=\tau_*(\pi, \phi, x)$ in (\ref{V3b}) such that $(\pi, \phi, x) \in C_*$.
        Here, the process $N^* = (N^*_t)_{t \ge 0}$ defined by:
        \begin{equation}
        \label{Mlocc4a}
        N^*_t = - \int_0^t \frac{\mu_1(X_s)-\mu_0(X_s)}{\sigma(X_s)} \, \frac{1-\pi_s}{\lambda + \alpha} \, \phi_s \, d{\bar B}_s
        \end{equation}
        is a continuous local martingale under $P_{\pi, \phi, x}$, and $\rho(x)$
        is the so-called {\it signal/noise ratio} function given by:
        \begin{equation}
        \label{rho4}
        \rho(x) = \left( \frac{\mu_1(x) - \mu_0(x)}{\sigma(x)} \right)^2
        \end{equation}
        for any $x > 0$. Observe that the assumption that the integral in (\ref{V3b}) taken up to the optimal stopping
        time $\tau_*$ is of finite expectation and the third inequality in (\ref{musig5}) yield that the integral in the left-hand
        side and the second integral on the right-hand side of (\ref{GIto4}) are of finite expectation too. Then, taking into
        account the fact that $0 \le \pi_t \le 1$ holds for all $t \ge 0$, and assuming that the process $(N^*_{\tau_* \wedge t})_{t \ge 0}$
        is a uniformly integrable martingale under $P_{\pi, \phi, x}$ (which is the case for the process $(M_{\tau_* \wedge t})_{t \ge 0}$
        from (\ref{tildeM3}) under the conditions of Lemma 3.3 below),
        by means of Doob's optional sampling theorem,
        we get from the expression in (\ref{GIto4}) that:
        \begin{equation}
        \label{GIto4a}
        E_{\pi, \phi, x} \int_0^{\tau_*} ({1-\pi_t}) \, \phi_t \, dt = E_{\pi, \phi, x} \int_0^{\tau_*}
        \frac{1-\pi_t}{\lambda + \alpha} \, d\phi_t
        - E_{\pi, \phi, x} \int_0^{\tau_*} \frac{1-\pi_t}{\lambda + \alpha} \, (\lambda + \rho(X_t) \, \pi_t \, \phi_t) \, dt
        \end{equation}
        is satisfied. Let us now take $x' > 0$ such that $x < x'$ and recall the fact that $(\pi_t, \phi_t, X_t)_{t \ge 0}$
        is a time-homogeneous Markov process. Assume that $(\pi, \phi, x) \in C_*$ is chosen sufficiently close
        to the stopping boundary $g_*(x)$, and note that $\tau_*=\tau_*(\pi, \phi, x)$ does not
        depend on $x'$. Hence, applying the comparison results from \cite{Ver} for solutions of stochastic differential equations,
        we obtain that the expression in (\ref{GIto4a}) yields that:
        \begin{align}
        \label{Vsigma4c''}
        V_*(\pi, \phi, x') 
        &\le E_{\pi, \phi, x'} \left[ 1 - \pi_{\tau_*}
        + \int_0^{\tau_*} (1-\pi_{t}) \, c \alpha \phi_{t} \, dt \right] \\
        \notag
        &\le E_{\pi, \phi, x} \left[ 1 - \pi_{\tau_*}
        + \int_0^{\tau_*} (1-\pi_{t}) \, c \alpha \phi_{t} \, dt \right] =
        V_*(\pi, \phi, x) 
        \end{align}
        holds, whenever $\rho(x)$ 
        is an increasing function on $(0, \infty)$.
        By virtue of the inequality in (\ref{Vsigma4c''}), we may therefore conclude that $(\pi, \phi, x') \in
        C_*$, so that the boundary $g_*(x)$ is increasing (decreasing) in (\ref{C5})-(\ref{D5}) whenever
        $\rho(x)$ is increasing (decreasing) on $(0, \infty)$, respectively.

        Summarizing the facts proved above, we are now ready to formulate the following assertion.

\vspace{6pt}

        {\bf Lemma 3.1.} {\it Suppose that $\mu_i(x)$, $i = 0, 1$, and $\sigma(x) > 0$
        are continuously differentiable functions on $(0, \infty)$ in (\ref{X3}) satisfying (\ref{musig5}).
        Assume that the process $(N^*_{\tau_* \wedge t})_{t \ge 0}$ from (\ref{Mlocc4a})
        is a uniformly integrable martingale. Then the optimal Bayesian time of alarm $\tau_*$
        in the quickest disorder detection problem of (\ref{V3b}) has the structure:
        \begin{equation}
        \label{tau4}
        \tau_* = \inf \{t \ge 0 \, | \, \phi_t \ge g_*(X_t) \}
        \end{equation}
        whenever the corresponding integral has finite expectation, so that $E_{\pi, \phi, x} \tau_* < \infty$ holds,
        for all $(\pi, \phi, x) \in [0, 1) \times [0, \infty) \times (0, \infty)$, and $\tau_* = 0$ otherwise.
        Moreover, 
        the property:
        \begin{equation}
        \label{C14}
        g_*(x): (0, \infty) \rightarrow (\lambda/(c \alpha), \infty) \;\,
        \text{is increasing/decreasing if} \;\, \rho(x) \;\, \text{is increasing/decreasing}
        \end{equation}
        holds with $\rho(x)$ defined in (\ref{rho4}), for all $x > 0$.}

\vspace{6pt}

        3.4. By means of standard arguments based on the application of
        It{\^o}'s formula, it is shown that the infinitesimal operator $\LL_{(\pi, \phi, X)}$
        of the process $(\pi_t, \phi_t, X_t)_{t \ge 0}$ from (\ref{pi3d}), (\ref{Phi2d}),
        and (\ref{X3a}) has the structure:
        \begin{align}
        \label{Lf3}
        &\LL_{(\pi, \phi, X)} = \lambda(1-\pi) \, \frac{\partial}{\partial \pi} +
        \bigg( \lambda + (\lambda + \alpha) \, \phi +
        \bigg( \frac{\mu_1(x)-\mu_0(x)}{\sigma(x)} \bigg)^2 \, \pi \, \phi \bigg) \, \frac{\partial}{\partial \phi}
        \\ \notag
        &\phantom{}
        + \big( \mu_0(x) + (\mu_1(x)-\mu_0(x)) \, \pi \big) \, \frac{\partial}{\partial x} +
        (\mu_1(x)-\mu_0(x)) \, \bigg( \pi(1-\pi) \, \frac{\partial^2}{\partial \pi \partial x}
        + \phi \, \frac{\partial^2}{\partial \phi \partial x} \bigg)\\ \notag
        &\phantom{}+ \frac{1}{2} \, \bigg( \frac{\mu_1(x)-\mu_0(x)}{\sigma(x)} \bigg)^2
        \bigg( \pi^2 (1-\pi)^2 \frac{\partial^2}{\partial \pi^2} +
        2 \pi(1-\pi) \, \phi \, \frac{\partial^2}{\partial \pi \partial \phi}
        + \phi^2 \, \frac{\partial^2}{\partial \phi^2} \bigg) + \frac{1}{2} \, {\sigma^2(x)} \, \frac{\partial^2}{\partial x^2}
        \end{align}
        for all $(\pi, \phi, x) \in [0, 1) \times [0, \infty) \times (0, \infty)$.

        According to the results of the general theory of optimal stopping problems for continuous
        time Markov processes (see, e.g. \cite{GS}, \cite[Chapter~III, Section~8]{S} and
        \cite[Chapter~IV, Section~8]{PSbook}), we can formulate the associated {\it free-boundary problem}
        for the unknown value function $V_*(\pi, \phi, x)$ from (\ref{V3b}) and the boundary $g_{*}(x)$
        from (\ref{tau4}):
        \begin{align}
        \label{VtpiC3}
        &(\LL_{(\pi, \phi, X)} V)(\pi, \phi, x) = - (1-\pi) \, c \alpha \phi  \quad
        \text{for} \quad (\pi, \phi, x) \in C \\
        \label{instop3}
        &V(\pi, \phi, x) \big|_{\phi=g(x)-} = 1 - \pi
        \quad \text{({\it instantaneous stopping})} \\
        \label{VGD3}
        &V(\pi, \phi, x) = 1 - \pi
        \quad \text{for} \quad (\pi, \phi, x) \in D \\
        \label{VGC3}
        &V(\pi, \phi, x) < 1 - \pi
        \quad \text{for} \quad (\pi, \phi, x) \in C
        \end{align}
        where $C$ and $D$ are defined as $C_*$ and $D_*$ in (\ref{C5}) and (\ref{D5})
        with $g(x)$ instead of $g_{*}(x)$, and the 
        condition in (\ref{instop3}) 
        is satisfied for all $\pi \in [0, 1)$ and $x > 0$.

        Note that the superharmonic characterization of the value function
        (see \cite{Dyn}, \cite[Chapter~III, Section~8]{S} and \cite[Chapter~IV, Section~9]{PSbook})
        implies that $V_*(\pi, \phi, x)$ from (\ref{V3b}) is the largest function
        satisfying (\ref{VtpiC3})-(\ref{VGC3}) with the boundary $g_{*}(x)$.

 \vspace{6pt}

      {\bf Remark 3.2.}
      Observe that, since the system in (\ref{VtpiC3})-(\ref{VGC3}) admits multiple solutions,
      we need to find some additional conditions which would specify the appropriate solution
      providing the value function and the optimal stopping boundary for the initial problem of (\ref{V3b}).
      In order to derive such conditions, we shall reduce the operator in (\ref{Lf3}) to the normal form.
      We also note that the fact that the stochastic differential equations for the posterior
      probability, the weighted likelihood ratio, and the observation process in (\ref{pi3d}), (\ref{Phi2d}), and
      (\ref{X3a}), respectively, are driven by the same (one-dimensional) innovation Brownian motion
      yields the property that the infinitesimal operator in (\ref{Lf3}) turns out to be of parabolic type.



       \vspace{6pt}

        3.5. In order to find the normal form of the operator in (\ref{Lf3})
        and formulate the appropriate optimal stopping and free-boundary
        problem, we use the one-to-one correspondence transformation of
        processes proposed by A.N. Kolmogorov in \cite{KolmII}. For this,
        let us define the process $Y = (Y_t)_{t \ge 0}$ by:
        \begin{equation}
        \label{Y4}
        Y_t = \log \phi_t
        - \int^{X_t}_{z} \frac{\mu_1(w)-\mu_0(w)}{\sigma^2(w)} \, dw
        \end{equation}
        for all $t \ge 0$, and any $z > 0$ fixed. Then, taking into account
        the assumption that the functions $\mu_i(x)$, $i = 0, 1$,
        and $\sigma(x)$ are continuously differentiable on $(0, \infty)$,
        by means of It{\^o}'s formula, we get that the process $Y$ admits the representation:
        \begin{equation}
        \label{dY3}
        dY_t = \Bigg( \frac{\lambda}{\phi_t} + \lambda + \alpha  - \frac{\sigma^2(X_t)}{2} \Bigg[ \frac{\mu^2_1(X_t)-\mu^2_0(X_t)}
        {\sigma^4(X_t)} + \frac{\partial}{\partial x} \left( \frac{\mu_1(x)-\mu_0(x)}{\sigma^2(x)} \right)
        \bigg|_{x = X_t} \Bigg] \Bigg) dt 
        \end{equation}
        with $Y_0 = y$ and
        \begin{equation}
        \label{y3}
        y = \log \phi
        - \int^{x}_{z} \frac{\mu_1(w)-\mu_0(w)}{\sigma^2(w)} \, dw
        \end{equation}
        for any $z > 0$ fixed.
        It is seen from the equation in (\ref{dY3}) that the
        process $Y$ started at $y \in \RR$ is of bounded variation.
%
        By virtue of the second inequality in (\ref{musig5}), 
        it follows from the relation in (\ref{Y4}) that there exists a one-to-one correspondence
        between the processes $(\pi, \phi, X)$ and $(\pi, \phi, Y)$. Hence, for any $z > 0$ fixed,
        the value function $V_*(\pi, \phi, x)$ from (\ref{V3b}) is equal to the one of the optimal
        stopping problem:
        \begin{equation}
        \label{V5c}
        U_*(\pi, \phi, y) = \inf_{\tau} E_{\pi, \phi, y}
        \left[ 1 - \pi_{\tau} + \int_0^{\tau} (1-\pi_t) \, c \alpha \phi_t \, dt \right]
        \end{equation}
        where the infimum is taken over all stopping times $\tau$ such that the integral is of finite expectation,
        so that $E_{\pi, \phi, y} \tau < \infty$ holds.
        Here, $P_{\pi, \phi, y}$ is a measure of the diffusion process $(\pi_t, \phi_t, Y_t)_{t \ge 0}$, started at some
        $(\pi, \phi, y) \in [0, 1) \times (0, \infty) \times \RR$ and solving the equations in (\ref{f2g}), (\ref{Phi2a}), and (\ref{Y4}).
        It thus follows from (\ref{C5})-(\ref{D5})
        that there exists a continuous function $h_{*}(y)$ such that
        $0 < \lambda/(c \alpha) \le h_*(y)$ holds for $y \in \RR$,
        and the optimal stopping time in the problem of (\ref{V5c}) has the structure:
        \begin{equation}
        \label{tau3y}
        \tau_* = \inf \{t \ge 0 \, | \, \phi_t \ge h_{*}(Y_t) \}
        \end{equation}
        whenever the corresponding integral is of finite expectation,
        so that $E_{\pi, \phi, y} \tau_* < \infty$ holds, and $\tau_* = 0$ otherwise.
        %


\vspace{6pt}

        3.6. Standard arguments then show that the infinitesimal operator $\LL_{(\pi, \phi, Y)}$
        of the process $(\pi, \phi, Y)$ from (\ref{pi3d}), (\ref{Phi2d}), and (\ref{dY3}) has the structure:
        \begin{align}
        \label{Lf3a}
        &\LL_{(\pi, \phi, Y)} = \lambda (1-\pi) \, \frac{\partial}{\partial \pi}
        + \bigg( \lambda + (\lambda + \alpha) \, \phi
        + \bigg( \frac{\mu_1(x(\phi, y))-\mu_0(x(\phi, y))}{\sigma(x(\phi, y))} \bigg)^2 \, \pi \, \phi \bigg)
        \, \frac{\partial}{\partial \phi} \\
        \notag
        &\phantom{}+ \frac{1}{2} \, \bigg( \frac{\mu_1(x(\phi, y))-\mu_0(x(\phi, y))}{\sigma(x(\phi, y))}
        \bigg)^2 \bigg( \pi^2 (1-\pi)^2 \frac{\partial^2}{\partial \pi^2} +
        2 \pi(1-\pi) \, \phi \, \frac{\partial^2}{\partial \pi \partial \phi}
        + \phi^2 \, \frac{\partial^2}{\partial \phi^2} \bigg) \\
        \notag
        &+ \left( \frac{\lambda}{\phi} + \lambda + \alpha - \frac{\sigma^2(x(\phi, y))}{2} \left[
        \frac{\mu^2_1(x(\phi, y))-\mu^2_0(x(\phi, y))}{\sigma^4(x(\phi, y))}
        + \frac{\partial}{\partial x}
        \left( \frac{\mu_1(x)-\mu_0(x)}{\sigma^2(x)} \right) \bigg|_{x = x(\phi, y)} \right] \right)
        \frac{\partial}{\partial y}
        \end{align}
        for all $(\pi, \phi, y) \in [0, 1) \times (0, \infty) \times \RR$. Here, because of the
        second inequality in (\ref{musig5}),
        the expression for $x(\phi, y) \equiv x(\phi, y; z)$ is uniquely determined
        by the relation in (\ref{y3}), for any $z > 0$.

        We are now ready to formulate the associated free-boundary problem
        for the unknown value function $U_*(\pi, \phi, y) \equiv U_*(\pi, \phi, y; z)$ from (\ref{V5c})
        and the boundary $h_{*}(y) \equiv h_{*}(y; z)$ from (\ref{tau3y}):
        \begin{align}
        \label{VtpiC5}
        &(\LL_{(\pi, \phi, Y)} U)(\pi, \phi, y) = - (1 - \pi) \, c \alpha \phi \quad
        \text{for} \quad \phi < h(y) \\
        \label{instop5}
        &U(\pi, \phi, y) \big|_{\phi=h(y)-} = 1 - \pi \quad (\text{\it instantaneous stopping}) \\
        \label{VGD5}
        &U(\pi, \phi, y) = 1 - \pi \quad \text{for} \quad \phi > h(y) \\
        \label{VGC5}
        &U(\pi, \phi, y) < 1 - \pi \quad \text{for} \quad \phi < h(y)
        \intertext{where
        the condition in (\ref{instop5}) 
        is satisfied for all $\pi \in [0, 1)$ and $y \in \RR$. 
        Moreover, we assume that the following conditions hold:}
        \label{smfit5}
        &\frac{\partial U}{\partial \phi}(\pi, \phi, y) \bigg|_{\phi=h(y)-} = 0
        \quad \text{({\it smooth fit})} \\
        \label{norent4}
        &\frac{\partial U}{\partial \phi}
        (\pi, \phi, y)\Big|_{\phi=0+} \quad \text{is finite}
        \intertext{and the one-sided derivative:}
        \label{party}
        &\frac{\partial U}{\partial y}(\pi, \phi, y) \bigg|_{\phi=h(y)-}  \quad \text{exists}
        \end{align}
        for all $\pi \in (0, 1)$, $y \in \RR$, and any $z > 0$ fixed.

        We further search for solutions of the parabolic-type free-boundary problem in (\ref{VtpiC5})-(\ref{VGC5})
        satisfying the conditions in (\ref{smfit5})-(\ref{party}) and such that the resulting boundaries
        are continuous and of bounded variation. Since such free-boundary problems cannot, in general,
        be solved explicitly, the existence and uniqueness of classical as well as viscosity solutions
        of the related variational inequalities and their connection with the optimal stopping problems
        have been extensively studied in the literature (see, e.g. \cite{Friedman}, \cite{BL}, \cite{Krylov} or \cite{O}).
        It particularly follows from the results of \cite[Chapter~XVI, Theorem~11.1]{Friedman} as well as 
        \cite[Chapter~V, Section~3, Theorem~14]{Krylov} with \cite[Chapter~VI, Section~4, Theorem~12]{Krylov} 
        that the free-boundary problem of (\ref{VtpiC5})-(\ref{VGC5}) with (\ref{smfit5})-(\ref{party}) admits a unique solution.



\vspace{6pt}

3.7. We continue with 
the following verification assertion related to the free-boundary
problem in (\ref{VtpiC5})-(\ref{party}).

\vspace{6pt}

     {\bf Lemma 3.3.} {\it
     Suppose that $\mu_i(x)$, $i = 0, 1$, and $\sigma(x) > 0$
     are continuously differentiable functions on $(0, \infty)$ in (\ref{X3})
     satisfying (\ref{musig5}).
     Assume that the function $U(\pi, \phi, y; h_*(y)) \equiv (1-\pi) H(\phi, y; h_*(y))$ and the
     continuous boundary of bounded variation $h_*(y)$ form a unique solution of the free-boundary problem
     in (\ref{VtpiC5})-(\ref{VGC5}) satisfying the conditions of (\ref{smfit5})-(\ref{party}).
     Then, the value function of the optimal stopping problem in (\ref{V5c}) takes the form:
     \begin{equation}
      \label{U*5}
      U_*(\pi, \phi, y)=
      \begin{cases}
      (1-\pi) H(\phi, y; h_*(y)), & \text{if} \; \; 0 \le \phi < h_*(y) \\
      1-\pi, & \text{if} \; \; \phi \ge h_*(y)
      \end{cases}
      \end{equation}
     and $h_*(y)$ provides the optimal stopping boundary for (\ref{tau3y}),
     whenever the corresponding integral is of finite expectation, so that
     $E_{\pi, \phi, y} \tau_* < \infty$ holds, for all
     $(\pi, \phi, y) \in [0, 1) \times (0, \infty) \times \RR$.}

        \vspace{6pt}

        {\bf Proof.} Let us denote by $U(\pi, \phi, y)$ the right-hand side of the expression in (\ref{U*5}).
        Hence, applying the change-of-variable formula with local time on surfaces from \cite{Pe1a} to
        $U(\pi, \phi, y)$ and $h_*(y)$, and taking into account the smooth-fit condition in (\ref{smfit5}),
        we obtain:
        \begin{equation}
        \label{rho3c}
        U(\pi_t, \phi_t, Y_t) = U(\pi, \phi, y) + \int_0^t (\LL_{(\pi, \phi, Y)} U) (\pi_s, \phi_s, Y_s) \, I(\phi_s \neq h_*(Y_s)) \, ds + M_t
        \end{equation}
        where the process $M = (M_t)_{t \ge 0}$ defined by:
        \begin{align}
        \label{tildeM3}
        M_t &= \int_0^t \frac{\partial U}{\partial \pi}(\pi_s, \phi_s, Y_s) \,
        \frac{\mu_1(X_s)-\mu_0(X_s)}{\sigma(X_s)} \, \pi_s (1-\pi_s) \, d{\bar B}_s \\
        \notag
        &\phantom{=\;\:} + \int_0^t \frac{\partial U}{\partial \phi}(\pi_s, \phi_s, Y_s)
        \, \frac{\mu_1(X_s)-\mu_0(X_s)}{\sigma(X_s)} \, \phi_s \, d{\bar B}_s
        \end{align}
        is a continuous local martingale under $P_{\pi, \phi, y}$ with respect to $({\cal F}_t)_{t \ge 0}$.


        It follows from the equation in (\ref{VtpiC5}) and the conditions of (\ref{VGD5})-(\ref{VGC5})
        that the inequality $(\LL_{(\pi, \phi, Y)} U)(\pi, \phi, y) \ge - (1-\pi) c \alpha \phi$
        holds for any $(\pi, \phi, y) \in [0, 1) \times (0, \infty) \times \RR$ such that
        $\phi \neq h_*(y)$, as well as $U(\pi, \phi, y) \le 1 - \pi$ is satisfied
        for all $(\pi, \phi, y) \in [0, 1) \times (0, \infty) \times \RR$.
        Recall the assumption that the boundary $h_*(y)$ is continuous and
        of bounded variation and the fact that the process $Y$ from (\ref{Y4}) is of bounded variation too.
        We thus conclude from the assumption of continuous differentiability of the functions $\mu_i(x)$, $i = 0, 1$,
        and $\sigma(x)$ that the time spent by the process $(\phi_t)_{t \ge 0}$ at the boundary $h_{*}(Y)$ is of Lebesgue
        measure zero, so that the indicator which appears in (\ref{rho3c})
        can be ignored.
        Hence, the expression in (\ref{rho3c}) yields that the inequalities:
        \begin{align}
        \label{rho3ea}
        &1-\pi_{\tau} + \int_0^{\tau} (1-\pi_t) \, c \alpha \phi_t \, dt \\
        \notag &\ge U(\pi_{\tau}, \phi_{\tau}, Y_{\tau})
        + \int_0^{\tau} (1-\pi_t) \, c \alpha \phi_t \, dt
        \ge U(\pi, \phi, y) + M_{\tau}
        \end{align}
        hold for any stopping time $\tau$ of the process $(\pi, \phi, Y)$
        started at $(\pi, \phi, y) \in [0, 1) \times (0, \infty) \times \RR$.

        Let $(\tau_n)_{n \in \NN}$ be an arbitrary localizing sequence of stopping times for the processes $M$.
        Taking the expectations with respect to the probability measure $P_{\pi, \phi, y}$ in (\ref{rho3ea}), by means of Doob's
        optional sampling theorem, we get that the inequalities:
        \begin{align}
        \label{rho3ec}
        &E_{\pi, \phi, y} \left[ 1-\pi_{\tau \wedge \tau_n} +
        \int_0^{\tau \wedge \tau_n} (1-\pi_t) \, c \alpha \phi_t \, dt \right] \\
        \notag
        &\ge E_{\pi, \phi, y} \left[ U(\pi_{\tau \wedge \tau_n}, \phi_{\tau \wedge \tau_n},
        Y_{\tau \wedge \tau_n}) + \int_0^{\tau \wedge \tau_n} (1-\pi_t) \, c \alpha \phi_t \, dt \right] \\
        \notag
        &\ge U(\pi, \phi, y) + E_{\pi, \phi, y} \big[ M_{\tau \wedge \tau_n} \big]
        = U(\pi, \phi, y)
        \end{align}
        hold for all $(\pi, \phi, y) \in [0, 1) \times (0, \infty) \times \RR$.
        Hence, letting $n$ go to infinity and using Fatou's lemma, we obtain:
        \begin{align}
        \label{rho3ee}
        &E_{\pi, \phi, y} \left[ 1-\pi_{\tau} +
        \int_0^{\tau} (1-\pi_t) \, c \alpha \phi_t \, dt \right] \\
        \notag
        &\ge E_{\pi, \phi, y} \left[ U(\pi_{\tau}, \phi_{\tau}, Y_{\tau})
        + \int_0^{\tau} (1-\pi_t) \, c \alpha \phi_t \, dt \right] \ge U(\pi, \phi, y)
        \end{align}
        for any stopping time $\tau$ and all $(\pi, \phi, y) \in [0, 1) \times (0, \infty) \times \RR$.
        By virtue of the structure of the stopping time in (\ref{tau3y}),
        it is readily seen that the inequalities in (\ref{rho3ee})
        hold with $\tau_*$ instead of $\tau$ when $\phi \ge h_*(y)$.

        It remains to show that the equalities are attained in (\ref{rho3ee})
        when $\tau_*$ replaces $\tau$, for $(\pi, \phi, y) \in [0, 1) \times
        (0, \infty) \times \RR$ such that $\phi < h_*(y)$.
        By virtue of the fact that the function $U(\pi, \phi, y)$ and the boundary $h_{*}(y)$
        satisfy the conditions in (\ref{VtpiC5}) and (\ref{instop5}), it follows from the expression
        in (\ref{rho3c}) and the structure of the stopping time in (\ref{tau3y}) that the equalities:
        \begin{equation}
        \label{rho3ga}
        U(\pi_{\tau_* \wedge \tau_n}, \phi_{\tau_* \wedge \tau_n}, Y_{\tau_* \wedge \tau_n})
        + \int_0^{\tau_* \wedge \tau_n} (1-\pi_t) \, c \alpha \phi_t \, dt
        = U(\pi, \phi, y) + M_{\tau_* \wedge \tau_n}
        \end{equation}
        hold for all $(\pi, \phi, y) \in [0, 1) \times (0, \infty) \times \RR$
        and any localizing sequence $(\tau_n)_{n \in \NN}$ of $M$.
        Hence, taking into account the assumption that the integral
        in (\ref{V3b}) taken up to $\tau_*$ is of finite expectation and
        using the fact that $0 \le U(\pi, \phi, y) \le 1$ holds,
        we conclude from the expression in (\ref{rho3ga}) that the process
        $(M_{\tau_* \wedge t})_{t \ge 0}$ is a uniformly integrable martingale.
        Therefore, taking the expectations in (\ref{rho3ga}) and letting $n$ go to infinity,
        we apply the Lebesgue dominated convergence theorem to obtain the equalities:
        \begin{align}
        \label{rho3ia}
        &E_{\pi, \phi, y}
        \left[ 1-\pi_{\tau_*} + \int_0^{\tau_*} (1-\pi_t) \, c \alpha \phi_t \, dt \right] \\
        \notag
        &= E_{\pi, \phi, y} \left[ U(\pi_{\tau_*}, \phi_{\tau_*}, Y_{\tau_*}) +
        \int_0^{\tau_*} (1-\pi_t) \, c \alpha \phi_t \, dt \right] = U(\pi, \phi, y)
        \end{align}
        for all $(\pi, \phi, y) \in [0, 1) \times (0, \infty) \times \RR$, which together with the inequalities
        in (\ref{rho3ee}) directly imply the desired assertion. $\square$

 \vspace{6pt}

        3.8. We are now in a position to formulate the main assertion of the
        paper, which follows from a straightforward combination of
        Lemma 3.3 above and standard change-of-variable
        arguments. More precisely, after obtaining the solution
        $U_*(\pi, \phi, y) \equiv (1-\pi) H_*(\phi, y; z)$
        with $h_{*}(y) \equiv h_{*}(y; z)$
        of the free-boundary problem in (\ref{VtpiC5})-(\ref{VGC5}),
        which satisfies the conditions in (\ref{smfit5})-(\ref{party}),
        we put $y = y(\pi, x; z)$ and $z = x$, in order to get the solution of the initial quickest
        detection problem with exponential penalty costs for a detection delay stated in (\ref{V3b}). 

 \vspace{6pt}

        {\bf Theorem 3.4.}
        {\it Suppose that the assumptions of Lemmas 3.1 and 3.3 hold.
        Then, in the quickest disorder detection problem of (\ref{V3b})
        for the observation process $X$ from (\ref{X3}), the Bayesian risk function takes the form
        $V_*(\pi, \phi, x) =
        U_*(\pi, \phi, y(\phi, x)) \equiv (1-\pi) H_*(\phi, y(\phi, x; x); x)$
        and the optimal stopping boundary $0 < \lambda/(c \alpha) \le g_{*}(x)$
        in (\ref{tau4}) satisfying (\ref{C14})
        is uniquely determined
        by the equation $g(x) = h_{*}(y(g(x), x)) \equiv h_{*}(y(g(x), x; x); x)$,
        for each $x > 0$ fixed. Here the function
        $U_*(\pi, \phi, y) \equiv (1-\pi) H_*(\phi, y; z)$ and the continuous
        boundary of bounded variation $h_{*}(y) \equiv h_{*}(y; z)$ form a unique
        solution of the free-boundary problem in (\ref{VtpiC5})-(\ref{party}),
        and the expression for $y(\phi, x) \equiv y(\phi, x; z)$
        is explicitly determined by the relation in (\ref{y3}),
        for all $(\pi, \phi, y) \in [0, 1) \times (0, \infty) \times \RR$ and any $z > 0$ fixed.}

\vspace{6pt}

        {\bf Remark 3.5.} 
        Observe that the optimal stopping time in the problem of (\ref{V3b})
        does not depend on the dynamics of the process $(\pi_t)_{t \ge 0}$, so that the
        two-dimensional process $(\phi_t, X_t)_{t \ge 0}$ turns out to be
        a sufficient statistic. 
        This fact is recognized as a consequence of the structure
        of the partial differential equation in (\ref{Lf3})-(\ref{VtpiC3}).
        However, the process $(\phi_t, X_t)_{t \ge 0}$ is not Markovian,
        and thus, in order to solve the optimal stopping problem of (\ref{V3b}),
        we need to add the component $(\pi_t)_{t \ge 0}$ and then
        operate with the resulting Markov process $(\pi_t, \phi_t, X_t)_{t \ge 0}$.

\vspace{6pt}

        Let us now give a short note concerning the case of bounded
        signal/noise ratio function $\rho(x)$ from (\ref{rho4}).

        \vspace{6pt}

        {\bf Remark 3.6.} Suppose that there exist some $0 < {\underline \rho} < {\overline \rho} < \infty$
        such that ${\underline \rho} \le \rho(x) \le {\overline \rho}$ holds for all $x > 0$.
        Let us denote by ${\underline V}_*(\pi, \phi, x)$ with ${\underline g}_*(x)$
        and by ${\overline V}_*(\pi, \phi, x)$ with ${\overline g}_*(x)$
        the solution of the Bayesian quickest disorder detection problem with exponential delay penalty,
        under $\rho(x) \equiv {\underline \rho}$ and $\rho(x) \equiv {\overline \rho}$, respectively.
        In those cases, the problem of (\ref{V3b}) degenerates into an optimal stopping problem
        for the two-dimensional Markov process $(\pi_t, \phi_t)_{t \ge 0}$, and the value functions
        ${\underline V}_*(\pi, \phi, x) \equiv {\underline V}_*(\pi, \phi)$ and ${\overline V}_*(\pi, \phi, x) \equiv {\overline V}_*(\pi, \phi)$
        with the stopping boundaries ${\underline g}_*(x) \equiv {\underline h}_*$ and ${\overline g}_*(x) \equiv {\overline h}_*$
        are given by the expressions in (\ref{Vg32}) and (\ref{Htil3a}) below, whenever $\eta=1/{\underline \rho}$ and
        $\eta=1/{\overline \rho}$, respectively.
        Taking into account the properties of the boundary $g_*(x)$ in (\ref{C14}) and the fact that
        $V_*(\pi, \phi, x) = 1 - \pi$ for all $\phi \ge g_*(x)$ and $0 \le \pi < 1$, we therefore conclude
        by standard comparison arguments that the inequalities ${\overline V}_*(\pi, \phi) \le V_*(\pi, \phi, x)
        \le {\underline V}_*(\pi, \phi)$ 
        and thus $0 < \lambda/(c \alpha) \le {\underline h}_* \le g_*(x) \le {\overline h}_*$
        hold for all $(\pi, \phi, x) \in [0, 1) \times [0, \infty) \times (0, \infty)$.


        \vspace{6pt}

        3.9. In order to pick up some special cases in which the
        free-boundary problem in (\ref{VtpiC5})-(\ref{party})
        can admit a simpler structure, for the rest of the section,
        in addition to the conditions in (\ref{musig5}), we suppose that the property:
        \begin{equation}
        \label{musigx3}
        \mu_i(x) = \frac{\eta_i \sigma^2(x)}{x} \quad \text{for some} \quad \eta_i \in \RR, \, i = 0, 1,
        \quad \text{such that} \quad \eta_0 \neq \eta_1 \quad \text{and} \quad \eta_0 + \eta_1 = 1
        \end{equation}
        holds for all $x > 0$.
        Moreover, we assume that the diffusion
        coefficient $\sigma(x)$ satisfies:
        \begin{equation}
        \label{musigas3}
        \sigma(x) \sim A_0 \, x^{\alpha} \quad \text{as} \quad x \downarrow 0 \qquad \text{and} \qquad
        \sigma(x) \sim A_{\infty} \, x^{\beta} \quad \text{as} \quad x \uparrow \infty
        \end{equation}
        with some $A_0, A_{\infty} > 0$ and $\alpha, \beta \in \RR$ such that
        $(1-\alpha) \eta \le 0$ and $(1-\beta) \eta \ge 0$ holds, where we set
        $\eta = 1/(\eta_1-\eta_0)$.
        In this case, the process $Y=(Y_t)_{t \ge 0}$ takes the form:
        \begin{equation}
        \label{Y5}
        Y_t = \log \phi_t
        - \frac{1}{\eta} \log \frac{X_t}{z} \equiv \log \phi + \int_0^t \bigg(\frac{\lambda}{\phi_s} + \lambda + \alpha \bigg) \, ds
        \quad \text{with} \quad \eta = \frac{1}{\eta_1-\eta_0}
        \end{equation}
        for any $z > 0$ fixed.
        It is easily seen from the structure of the expression in (\ref{Y5}) that the one-to-one
        correspondence between the processes $(\pi_t, \phi_t, X_t)_{t \ge 0}$ and $(\pi_t, \phi_t, Y_t)_{t \ge 0}$ remains true in this case.
        Hence, getting the expression for $X_t$ from (\ref{Y5}) and substituting it
        into the equations of (\ref{pi3d}) and (\ref{Phi2d}), we obtain:
        \begin{equation}
        \label{pi4d}
        d\pi_t = \lambda (1 - \pi_t) \, dt + \frac{\sigma (z e^{- \eta Y_t} \phi_t^{\eta})}
        {\eta z e^{-\eta Y_t} \phi_t^{\eta}} \, \pi_t (1 - \pi_t) \, d{\bar B}_t 
        \end{equation}
        with $\pi_0=\pi$ and
        \begin{equation}
        \label{phi4d}
        d\phi_t = \bigg( \lambda + (\lambda + \alpha) \phi_t + \frac{\sigma^2 (z e^{- \eta Y_t} \phi_t^{\eta})}
        {\eta^2 z^2 e^{- 2 \eta Y_t} \phi_t^{2 \eta}} \, \pi_t \, \phi_t \bigg) \, dt + \frac{\sigma (z e^{- \eta Y_t} \phi_t^{\eta})}
        {\eta z e^{-\eta Y_t} \phi_t^{\eta}} \, \phi_t \, d{\bar B}_t 
        \end{equation}
        with $\phi_0=\phi$, for any $z > 0$ fixed.
        Applying It{\^o}'s formula to the expression in (\ref{Y5})
        and taking into account the representations in (\ref{pi3d}) and (\ref{X3a})
        as well as the assumption of (\ref{musigx3}), we get:
        \begin{equation}
        \label{Y3b}
        dY_t = \bigg(\frac{\lambda}{\phi_t} + \lambda + \alpha \bigg) \, dt
        \end{equation}
        with $Y_0 = y$. It thus follows that the infinitesimal
        operator $\LL_{(\pi, \phi, Y)}$ from (\ref{Lf3a}) takes the form:
        \begin{align}
        \label{Lf5}
        &\LL_{(\pi, \phi, Y)} = \lambda (1-\pi) \, \frac{\partial}{\partial \pi}
        + \bigg( \lambda + (\lambda + \alpha) \phi
        +  \frac{\sigma^2 (z e^{- \eta y} \phi^{\eta})}{\eta^2 z^2 e^{- 2 \eta y} \phi^{2 \eta}}
        \, \pi \, \phi \bigg) \, \frac{\partial}{\partial \phi} \\ \notag
        &\phantom{}+ \frac{1}{2} \, \frac{\sigma^2 (z e^{- \eta y} \phi^{\eta})}
        {\eta^2 z^2 e^{- 2 \eta y} \phi^{2 \eta}} \, \bigg( \pi^2 (1-\pi)^2 \frac{\partial^2}{\partial \pi^2} +
        2 \pi(1-\pi) \, \phi \, \frac{\partial^2}{\partial \pi \partial \phi}
        + \phi^2 \, \frac{\partial^2}{\partial \phi^2} \bigg) + \bigg(\frac{\lambda}{\phi} + \lambda + \alpha \bigg)
        \, \frac{\partial}{\partial y}
        \end{align}
        for all $(\pi, \phi, y) \in [0, 1) \times (0, \infty) \times \RR$ and any $z > 0$ fixed.




\vspace{6pt}

        3.10. Let us now introduce the function ${\hat U}(\pi, \phi, y) \equiv (1-\pi) {\hat H}(\phi, y)$
        and the boundary ${\hat h}(y)$ as a solution of the free-boundary problem consisting
        of the differential equation:
        \begin{equation}
        \label{hatVtpiC5}
        \bigg( \big( \lambda + (\lambda + \alpha) \phi \big) \, \frac{\partial H}{\partial \phi}
        + \frac{1}{2} \, \frac{\sigma^2 (z e^{- \eta y} \phi^{\eta})}{\eta^2 z^2 e^{- 2 \eta y} \phi^{2 \eta}} \,
        \phi^2 \, \frac{\partial^2 H}{\partial \phi^2} - \lambda H \bigg) (\phi, y) = - c \alpha \phi
        \quad \text{for} \quad \phi < h(y)
        \end{equation}
        instead of the one in (\ref{VtpiC5}), for each $y > 0$ fixed, and the
        conditions of (\ref{instop5})-(\ref{VGC5}) as well as (\ref{smfit5})-(\ref{party}).
        The general solution of the resulting second-order {\it ordinary} differential equation in (\ref{hatVtpiC5}) takes the form:
        \begin{equation}
        \label{genH}
        H(\phi, y) = C_0(y) \, H_0(\phi, y) + C_{\infty}(y) \, H_{\infty}(\phi, y) - c (1 + \phi)
        \end{equation}
        where
        $H_i(\phi, y)$, $i = 0, \infty$,
        form a system of fundamental positive
        solutions (i.e. nontrivial linearly independent particular solutions) of the corresponding
        {\it homogeneous} differential equation, and $C_i(y)$, $i = 0, \infty$, are some arbitrary continuously
        differentiable functions, so that the condition in (\ref{party}) holds. By virtue of the
        assumptions of (\ref{musig5}) and 
        taking into account the arguments from \cite[Section~4]{GS1}, we can identify by
        $H_0(\phi, y)$ a decreasing solution that has a singularity at zero and by
        $H_{\infty}(\phi, y)$ an increasing solution that has a singularity at infinity.

        Observe that we should have $C_0(y) = 0$ in (\ref{genH}),
        since otherwise $U(\pi, \phi, y) \equiv (1-\pi) H(\phi, y) \to \pm \infty$ as $\phi \downarrow 0$,
        that must be excluded by virtue of the obvious fact that the value function in (\ref{V5c})
        is bounded at $\phi = 0$, for any $y \in \RR$ fixed. Then, applying the conditions
        of (\ref{instop5}) and (\ref{smfit5}) to the function in (\ref{genH}) with $C_0(y) = 0$,
        we get that the equalities:
        \begin{equation}
        \label{genHa}
        C_{\infty}(y) \, H_{\infty}(h(y), y) = c (1 + h(y)) + 1 \quad \text{and} \quad
        C_{\infty}(y) \, \frac{\partial H_{\infty}}{\partial \phi}(\phi, y)
        \bigg|_{\phi = h(y)} = c
        \end{equation}
        hold for $y \in \RR$ fixed. Hence, solving the equations of (\ref{genHa}),
        we get that the solution of the system of (\ref{hatVtpiC5})
        with (\ref{instop5}) and (\ref{smfit5})-(\ref{norent4}) is given by:
        \begin{equation}
        \label{LH3b}
        H(\phi, y; {\hat h}(y)) = \big(c (1 + {\hat h}(y)) + 1 \big) \,
        \frac{H_{\infty}(\phi, y)}{H_{\infty}({\hat h}(y), y)} - c (1 + \phi)
        \end{equation}
        for all $0 \le \phi < {\hat h}(y)$, so that $0 \le H(\phi, y; {\hat h}(y))
        \equiv H(\phi, y; z; {\hat h}(y; z)) \le 1$ holds, where ${\hat h}(y)$ satisfies the equation:
        \begin{equation}
        \label{Htil3a}
        \frac{\partial H_{\infty}}{\partial \phi}(\phi, y) \bigg|_{\phi = h(y)}
        = \frac{c H_{\infty}(h(y), y)}{c(1+h(y))+1}
        \end{equation}
        for any $y \in \RR$ fixed.


        Taking into account the facts proved above, let us formulate the
        following assertion.

        \vspace{6pt}

        {\bf Corollary 3.7.}
        {\it Suppose that $\mu_i(x)$, $i = 0, 1$, and $\sigma(x) > 0$ are continuously differentiable functions
        on $(0, \infty)$ in (\ref{X3}) satisfying (\ref{musig5}) and (\ref{musigx3})-(\ref{musigas3})
        with $\alpha, \beta \in \RR$ such that $(1-\alpha) \eta \le 0$ and $(1-\beta) \eta \ge 0$,
        where $\eta = 1/(\eta_1-\eta_0)$.
        Assume that ${\hat h}(y)$ provides a unique solution of the equation in (\ref{Htil3a}) for all $y \in \RR$.
        Then, using the same arguments as in the proof of Lemma 3.3 above, it is shown that the function:
        \begin{equation}
        \label{Vg32}
        {\hat U}(\pi, \phi, y) \equiv (1-\pi){\hat H}(\phi, y) \quad \text{with} \quad
        {\hat H}(\phi, y) =
        \begin{cases}
        H(\phi, y; {\hat h}(y)), & \text{if} \; \; 0 \le \phi < {\hat h}(y) \\
        1, & \text{if} \; \; \phi \ge {\hat h}(y)
        \end{cases}
        \end{equation}
        coincides with the value function of the optimal stopping problem:
        \begin{align}
        \label{hatV5c}
        &{\hat U}(\pi, \phi, y) \\
        \notag
        &= \inf_{\tau} E_{\pi, \phi, y} \bigg[ 1 - \pi_{\tau} +
        \int_0^{\tau} (1-\pi_t) \, \bigg( c \alpha \phi_t - \bigg(\frac{\lambda}{\phi_t} + \lambda + \alpha \bigg) \,
        \frac{\partial {\hat H}}{\partial y}(\phi_t, Y_t) \, I(\phi_t < {\hat h}(Y_t)) \bigg) \, dt \bigg]
        \end{align}
        which corresponds to the Bayesian risk function in (\ref{V5c}).
        Moreover, ${\hat h}(y) \equiv {\hat h}(y; z)$ determined by (\ref{Htil3a})
        provides a hitting boundary for the stopping time:
        \begin{equation}
        \label{hattau3}
        {\hat \tau} = \inf \{t \ge 0 \, | \, \phi_t \ge {\hat h}(Y_t) \}
        \end{equation}
        which turns out to be optimal in (\ref{hatV5c}) whenever the
        integral above is of finite expectation, and ${\hat \tau} = 0$ otherwise,
        for any $z > 0$ fixed.}

        \vspace{6pt}

        {\bf Remark 3.8.} Note that the function ${\hat U}(\pi, \phi, y)$ in (\ref{hatV5c})
        and the boundary ${\hat h}(y)$ in (\ref{hattau3}) provide lower (upper) and upper
        (lower) estimates for the initial value function $U_*(\pi, \phi, y)$ in (\ref{V5c})
        and the optimal stopping boundary $h_*(y)$ in (\ref{tau3y}), whenever the function
        $y \mapsto {\hat H}(\phi, y)$ is increasing (decreasing) on $\RR$.
        According to Remark 3.6 above and the structure of the change of variables in (\ref{y3}),
        such a situation occurs when $\rho(x)$ from (\ref{rho4}) is an increasing (decreasing)
        function on $(0, \infty)$ and $\eta_0 < \eta_1$ ($\eta_0 > \eta_1$) in (\ref{musigx3}),
        respectively.

\section{\dot The case of linear delay penalty costs}

        In this section, we provide some results, which are related to the quickest
        detection problem with linear delay penalty costs of Example 2.2 above.

        \vspace{6pt}

        4.1. Following the arguments of Subsection 3.1 above and applying Doob's optional sampling theorem,
        we get from (\ref{GIto3}) that the equality:
        \begin{equation}
        \label{2.13'}
        E_{\pi, \varphi, x} \bigg[ 1 - \pi_{\tau} + \int_0^{\tau} (1 - \pi_t) \, c \varphi_t \, dt \bigg]
        = 1 - \pi + E_{\pi, \varphi, x} \int_0^{\tau} (1 - \pi_t) \, (c \varphi_t - \lambda) \, dt
        \end{equation}
        holds for all $(\pi, \varphi, x) \in [0, 1) \times [0, \infty) \times (0, \infty)$
        and any stopping time $\tau$ satisfying $E_{\pi, \varphi, x} < \infty$.
        Taking into account the structure of the reward in (\ref{V30}),
        it is also seen from (\ref{2.13'}) that it is never optimal
        to stop when $\varphi_t < \lambda/c$ for any $t \ge 0$.
        This shows that all the points $(\pi, \varphi, x)$ such that
        $\varphi < \lambda/c$ belong to the continuation region:
        \begin{equation}
        \label{C3aa'}
        C' = \{ (\pi, \varphi, x) \in [0, 1) \times [0, \infty) \times (0, \infty) \, | \, V'(\pi, \varphi, x) <  1 - \pi \}.
        \end{equation}
        Then, combining the arguments in \cite[Chapter~IV, Section~3]{S} with the ones in Subsection 3.2 above,
        we obtain that the continuation region in (\ref{C3aa'}) for the optimal stopping problem of (\ref{V30}) takes the form:
        \begin{equation}
        \label{C5'}
        C' = \{(\pi, \varphi, x) \in [0, 1) \times [0, \infty) \times (0, \infty) \, | \, \varphi < g'(x) \}
        \end{equation}
        so that the corresponding stopping region is the closure of the set:
        \begin{equation}
        \label{D5'}
        D' = \{(\pi, \varphi, x) \in [0, 1) \times [0, \infty) \times (0, \infty) \, | \, \varphi > g'(x) \}.
        \end{equation}

        \vspace{3pt}

        4.2. In order to characterize the behavior of the boundary $g'(x)$ in (\ref{C5'})-(\ref{D5'}),
        we observe from the equation in (\ref{phi3d}) that the expression:
        \begin{equation}
        \label{GIto5}
        \int_0^{\tau'} ({1-\pi_t}) \, \varphi_t \, dt
        = \int_0^{\tau'} \frac{1-\pi_t}{\lambda} \, d\varphi_t
        - \int_0^{\tau'} \frac{1-\pi_t}{\lambda} \, (\lambda + \rho(X_t) \, \pi_t \, \varphi_t) \, dt + N_{\tau'}'
        \end{equation}
        holds for the optimal stopping time $\tau'=\tau'(\pi, \varphi, x)$ in (\ref{V30}) such that $(\pi, \varphi, x) \in C'$.
        Here the process $N' = (N_t')_{t \ge 0}$ defined by:
        \begin{equation}
        \label{Mlocc5a}
        N_t' = - \int_0^t \frac{\mu_1(X_s)-\mu_0(X_s)}{\sigma(X_s)} \, \frac{1-\pi_s}{\lambda} \, \varphi_s \, d{\bar B}_s
        \end{equation}
        is a continuous local martingale under $P_{\pi, \varphi, x}$, and
        the function $\rho(x)$ is given by (\ref{rho4}). Note that the assumption that
        $E_{\pi, \varphi, x} \tau' < \infty$ holds and the third inequality in (\ref{musig5}) yield that
        the integral in the left-hand side and the second integral on the right-hand side of (\ref{GIto5})
        are of finite expectation. Moreover, by virtue of the facts that $(1-\pi_t) \varphi_t = \pi_t$
        and $0 \le \pi_t \le 1$ holds for all $t \ge 0$, and taking into account the third inequality
        in (\ref{musig5}), we see from (\ref{Mlocc5a}) that the process $(N_{\tau' \wedge t}')_{t \ge 0}$
        is a uniformly integrable martingale under $P_{\pi, \varphi, x}$.
        Then, applying Doob's optional sampling theorem, we get from the expression in (\ref{GIto5}) that:
        \begin{equation}
        \label{GIto5a}
        E_{\pi, \varphi, x} \int_0^{\tau'} (1-\pi_t) \, \varphi_t \, dt = E_{\pi, \varphi, x} \int_0^{\tau'}
        \frac{1-\pi_t}{\lambda} \, d\varphi_t - E_{\pi, \varphi, x} \int_0^{\tau'}
        \frac{1-\pi_t}{\lambda} \, (\lambda + \rho(X_t) \, \pi_t \, \varphi_t) \, dt
        \end{equation}
        is satisfied.
        Let us now take $x' > 0$ such that $x < x'$ and recall the fact that $(\pi_t, \varphi_t, X_t)_{t \ge 0}$
        is a time-homogeneous Markov process. Assume that $(\pi, \varphi, x) \in C'$ is chosen sufficiently close
        to the stopping boundary $g'(x)$, and note that $\tau'=\tau'(\pi, \varphi, x)$ does not depend on $x'$.
        Hence, by means of the comparison results for solutions of stochastic differential equations,
        we obtain that the expression in (\ref{GIto5a}) yields:
        \begin{align}
        \label{Vsigma5c''}
        V'(\pi, \varphi, x') - (1 - \pi) &\le E_{\pi, \varphi, x'} \left[ 1 - \pi_{\tau_*}
        + \int_0^{\tau_*} (1-\pi_{t}) \, c \varphi_{t} \, dt \right] \\
        \notag
        &\le E_{\pi, \varphi, x} \left[ 1 - \pi_{\tau_*}
        + \int_0^{\tau_*} (1-\pi_{t}) \, c \varphi_{t} \, dt \right] =
        V'(\pi, \varphi, x) - (1-\pi)
        \end{align}
        whenever $\rho(x)$
        is an increasing function on $(0, \infty)$.
        By virtue of the inequality in (\ref{Vsigma5c''}), we may therefore conclude that $(\pi, \varphi, x') \in C'$,
        so that the boundary $g'(x)$ is increasing (decreasing) in (\ref{C5'})-(\ref{D5'}) whenever $\rho(x)$
        is increasing (decreasing) on $(0, \infty)$, respectively.


        Summarizing the facts proved above, we now formulate the assertions related to the Bayesian quickest detection
        problem with linear penalty costs for a detection delay, which are proved using the arguments from the previous
        section.

        \vspace{6pt}

        {\bf Lemma 4.1.} {\it Suppose that $\mu_i(x)$, $i = 0, 1$, and $\sigma(x) > 0$
        are continuously differentiable functions on $(0, \infty)$ in (\ref{X3}) satisfying (\ref{musig5}).
        Then the optimal Bayesian time of alarm $\tau'$
        in the quickest disorder detection problem (\ref{V30}) has the structure:
        \begin{equation}
        \label{tau5'}
        \tau' = \inf \{ t \ge 0 \, | \, \varphi_t \ge g'(X_t) \}
        \end{equation}
        whenever $E_{\pi, \varphi, x} \tau' < \infty$ holds, for all
        $(\pi, \varphi, x) \in [0, 1) \times [0, \infty) \times (0, \infty)$, and $\tau' = 0$ otherwise.
        Moreover, 
        the property:
        \begin{equation}
        \label{C15}
        g'(x): (0, \infty) \rightarrow (\lambda/c, \infty) \;\,
        \text{is increasing/decreasing if} \;\, \rho(x) \;\, \text{is increasing/decreasing}
        \end{equation}
        holds with $\rho(x)$ defined in (\ref{rho4}), for all $x > 0$.}

        \vspace{6pt}


        \vspace{6pt}

        {\bf Theorem 4.2.}
        {\it Suppose that the assumptions of Lemmas 4.1 and 3.3 hold with $\phi = \varphi$, $\alpha = 0$ in
        (\ref{Lf3a}), and $\alpha = 1$ in (\ref{VtpiC5}).
        Then, in the quickest disorder detection problem of (\ref{V30})
        for the observation process $X$ from (\ref{X3}), the Bayesian risk function takes the form
        $V'(\pi, \varphi, x) =
        U'(\pi, \varphi, y(\varphi, x)) \equiv (1-\pi) H'(\varphi, y(\varphi, x; x); x)$
        and the optimal stopping boundary $0 < \lambda/c \le g'(x)$
        in (\ref{tau5'}) satisfying (\ref{C15})
        is uniquely determined
        by the equation $g(x) = h'(y(g(x), x)) \equiv h'(y(g(x), x; x); x)$,
        for each $x > 0$ fixed. Here the function
        $U'(\pi, \varphi, y) \equiv (1-\pi) H'(\varphi, y; z)$ and the bounded continuous boundary
        of bounded variation $h'(y) \equiv h'(y; z)$
        form a unique solution of the free-boundary problem in (\ref{VtpiC5})-(\ref{party}),
        and the expression for $y(\varphi, x) \equiv y(\varphi, x; z)$
        is explicitly determined by the relation in (\ref{y3}) with $\phi = \varphi$,
        for all $(\pi, \varphi, y) \in [0, 1) \times (0, \infty) \times \RR$ and any $z > 0$ fixed.}



        \vspace{6pt}

        {\bf Remark 4.3.} Suppose that there exist some $0 < {\underline \rho} < {\overline \rho} < \infty$
        such that ${\underline \rho} \le \rho(x) \le {\overline \rho}$ holds for all $x > 0$.
        Let us denote by ${\underline V}'(\pi, \varphi, x)$ with ${\underline g}'(x)$
        and by ${\overline V}'(\pi, \varphi, x)$ with ${\overline g}'(x)$
        the solution of the Bayesian quickest disorder detection problem with linear delay penalty,
        under $\rho(x) \equiv {\underline \rho}$ and $\rho(x) \equiv {\overline \rho}$, respectively.
        In those cases, the problem of (\ref{V30}) degenerates into an optimal stopping problem
        for the one-dimensional Markov process $(\pi_t)_{t \ge 0}$ being equivalent to $(\varphi_t)_{t \ge 0}$, and the value
        functions ${\underline V}'(\pi, \varphi, x) \equiv {\underline V}'(\pi, \varphi) \equiv {\underline V}'(\varphi/(1+\varphi),
        \varphi)$ and ${\overline V}'(\pi, \varphi, x) \equiv {\overline V}'(\pi, \varphi) \equiv {\overline V}'(\varphi/(1+\varphi),
        \varphi)$ with the stopping boundaries ${\underline g}'(x) \equiv {\underline h}'$ and
        ${\overline g}'(x) \equiv {\overline h}'$ are given by the expressions in (\ref{Vh32}) and (\ref{Htil3h}) below,
        whenever $\eta=1/{\underline \rho}$ and $\eta=1/{\overline \rho}$, respectively.
        Taking into account the properties of the boundary $g'(x)$ in (\ref{C15}) and the fact that
        $V'(\pi, \varphi, x) = 1 - \pi$ for all $\varphi \ge g'(x)$ and $0 \le \pi < 1$, we therefore conclude
        by standard comparison arguments that the inequalities ${\overline V}'(\pi, \varphi) \le V'(\pi, \varphi, x)
        \le {\underline V}'(\pi, \varphi)$ 
        and thus $0 < \lambda/c \le {\underline h}' \le g'(x) \le {\overline h}'$
        hold for all $(\pi, \varphi, x) \in [0, 1) \times [0, \infty) \times (0, \infty)$.


        \vspace{6pt}

        4.4. Let us finally introduce the function ${\tilde U}(\varphi/(1+\varphi), \varphi, y) \equiv {\tilde G}(\varphi, y)$
        and the boundary ${\hat h}(y)$ as a solution of the free-boundary problem consisting
        of the differential equation:
        \begin{equation}
        \label{hatVtpiC5a}
        \bigg( \lambda (1 + \varphi) \, \frac{\partial G}{\partial \varphi} +
        \frac{\sigma^2 (z e^{- \eta y} \varphi^{\eta})}{\eta^2 z^2 e^{- 2 \eta y} \varphi^{2 \eta}}
        \bigg(\frac{\varphi^2}{1+\varphi} \, \frac{\partial G}{\partial \varphi} + \frac{\varphi^2}{2}
        \, \frac{\partial^2 G}{\partial \varphi^2} \bigg) \bigg) (\varphi, y) = - \frac{c \varphi}{1+\varphi}
        \quad \text{for} \quad \varphi < h(y)
        \end{equation}
        instead of the one in (\ref{VtpiC5}), for each $y > 0$ fixed, and the
        conditions of (\ref{instop5})-(\ref{VGC5}) as well as (\ref{smfit5})-(\ref{party})
        with $\phi = \varphi$ and $\pi = \varphi/(1+\varphi)$.
        The general solution of the resulting first-order linear ordinary differential equation
        for $\varphi \mapsto ({\partial G}/{\partial \varphi})(\varphi, y)$ takes the form:
        \begin{align}
        \label{hatVtpiC5b}
        \frac{\partial G}{\partial \varphi}(\varphi, y) &= \frac{C(y)}{(1+\varphi)^2} \,
        \exp \left( \int_{\varphi}^w \frac{\lambda(1+u)}{u^2} \, \frac{2 \eta^2 z^2 e^{- 2 \eta y} u^{2 \eta}}
        {\sigma^2(z e^{-\eta y} u^{\eta})} \, du \right) \\
        \notag
        &\phantom{=\:\;} - \int_0^{\varphi} \frac{c (1+u)}{u (1+\varphi)^2} \, \frac{2 \eta^2 z^2 e^{- 2 \eta y}
        u^{2 \eta}} {\sigma^2(z e^{-\eta y} u^{\eta})} \, \exp \left( - \int_u^{\varphi}
        \frac{\lambda(1+v)}{v^2} \, \frac{2 \eta^2 z^2 e^{- 2 \eta y} v^{2 \eta}}
        {\sigma^2(z e^{-\eta y} v^{\eta})} \, dv \right) du
        \end{align}
        where $C(y)$ is an arbitrary continuously differentiable function, for each $y \in \RR$ and any $z, w > 0$ fixed.
        By virtue of the assumptions of (\ref{musig5}), 
        we see that the term in the first line of (\ref{hatVtpiC5b}) above tends to infinity as $\varphi \downarrow 0$,
        so that $({\partial G}/{\partial \varphi})(\varphi, y) \to \pm \infty$ as $C(y) \neq 0$, for any $y \in \RR$ fixed.
        We should thus choose $C(y) = 0$, that is equivalent to the property in (\ref{norent4}).
        Hence, integrating the equation in (\ref{hatVtpiC5b}), we therefore obtain that the solution of the system
        of (\ref{hatVtpiC5a}) with (\ref{instop5}) and (\ref{smfit5})-(\ref{norent4}) is given by:
        \begin{align}
        \label{LH3h}
        &G(\varphi, y; {\tilde h}(y)) = 1/(1+{\tilde h}(y)) \\
        \notag
        &+ \int_{\varphi}^{{\tilde h}(y)} \int_0^{w}
        \frac{c(1+u)}{u(1+w)^2} \, \frac{2 \eta^2 z^2 e^{- 2 \eta y} u^{2 \eta}}
        {\sigma^2(z e^{-\eta y} u^{\eta})} \, \exp \left( - \int_u^{w}
        \frac{\lambda(1+v)}{v^2} \, \frac{2 \eta^2 z^2 e^{- 2 \eta y} v^{2 \eta}}
        {\sigma^2(z e^{-\eta y} v^{\eta})} \, dv \right) du \, dw
        \end{align}
        for all $0 \le \varphi < {\tilde h}(y)$, so that
        $0 \le G(\varphi, y; {\tilde h}(y)) \equiv G(\varphi, y; z; {\tilde h}(y; z)) \le 1/(1+\varphi)$ holds,
        where ${\tilde h}(y)$ satisfies the equation:
        \begin{equation}
        \label{Htil3h}
        \int_0^{h(y)} \frac{c(1+u)}{u} \, \frac{2 \eta^2 z^2 e^{- 2 \eta y} u^{2 \eta}}
        {\sigma^2(z e^{-\eta y} u^{\eta})} \, \exp \bigg( - \int_u^{h(y)}
        \frac{\lambda(1+v)}{v^2} \, \frac{2 \eta^2 z^2 e^{- 2 \eta y} v^{2 \eta}}
        {\sigma^2(z e^{-\eta y} v^{\eta})} \, dv \bigg) du = 1
        \end{equation}
        for each $y \in \RR$ and any $z > 0$ fixed.



        Summarizing these facts above, let us formulate the following assertion.

        \vspace{6pt}

        {\bf Corollary 4.4.}
        {\it Suppose that $\mu_i(x)$, $i = 0, 1$, and $\sigma(x) > 0$ are continuously differentiable functions
        on $(0, \infty)$ in (\ref{X3}) satisfying (\ref{musig5}) and (\ref{musigx3})-(\ref{musigas3}) with
        $\alpha, \beta \in \RR$ such that $(1-\alpha) \eta \le 0$ and $(1-\beta) \eta \ge 0$, where $\eta = 1/(\eta_1-\eta_0)$.
        Assume that ${\tilde h}(y)$ provides a unique solution of the equation in (\ref{Htil3h}) for all $y \in \RR$.
        Then, using the same arguments as in the proof of Lemma 3.3 above, it is shown that
        the function: 
        \begin{equation}
        \label{Vh32}
        {\tilde U}(\pi, \varphi, y) \equiv {\tilde G}(\varphi, y) =
        \begin{cases}
        G(\varphi, y; {\tilde h}(y)), & \text{if} \; \; 0 \le \varphi < {\tilde h}(y) \\
        1/(1+\varphi), & \text{if} \; \; \varphi \ge {\tilde h}(y)
        \end{cases}
        \end{equation}
        coincides with the value function of the optimal stopping problem:
        \begin{equation}
        \label{hatV5h}
        {\tilde U}(\pi, \varphi, y)
        = \inf_{\tau} E_{\pi, \varphi, y} \bigg[ \frac{1}{1+\varphi_{\tau}}
        + \int_0^{\tau} 
        \bigg( \frac{c \varphi_t}{1+\varphi_t}
        - \bigg( \frac{\lambda}{\varphi_t} + \lambda \bigg)
        \frac{\partial {\tilde G}}{\partial y}(\varphi_t, Y_t)
        \, I({\varphi}_t < {\tilde h}(Y_t)) \bigg) \, dt \bigg]
        \end{equation}
        with $\pi = \varphi/(1+\varphi)$, which corresponds to the Bayesian risk function in (\ref{V30}).
        Moreover, $0 < \lambda/c \le {\tilde h}(y) \equiv {\tilde h}(y; z)$
        determined by (\ref{Htil3h}) provides a hitting boundary for the stopping time:
        \begin{equation}
        \label{hattau3h}
        {\tilde \tau} = \inf \{t \ge 0 \, | \, \varphi_t \ge {\tilde h}(Y_t) \}
        \end{equation}
        which turns out to be optimal in (\ref{hatV5h})
        whenever the integral above is of finite expectation,
        for any $z > 0$ fixed.}

        \vspace{6pt}

        {\bf Remark 4.5.} Note that the function ${\tilde U}(\pi, \varphi, y)$ in (\ref{hatV5h})
        and the boundary ${\tilde h}(y)$ in (\ref{hattau3h}) provide lower (upper) and upper
        (lower) estimates for the initial value function $U'(\pi, \varphi, y)$ defined as in
        (\ref{V5c}) with $\alpha = 1$ and the optimal stopping boundary $h'(y)$ defined as in
        (\ref{tau3y}) with $\phi = \varphi$,
        whenever the function $y \mapsto {\tilde G}(\varphi, y)$ is increasing (decreasing) on $\RR$.
        According to Remark 4.1 and the structure of the change of variables in (\ref{y3}) with $\phi = \varphi$,
        such a situation occurs when $\rho(x)$ from (\ref{rho4}) is an increasing (decreasing)
        function on $(0, \infty)$ and $\eta_0 < \eta_1$ ($\eta_0 > \eta_1$) in (\ref{musigx3}),
        respectively.

\vspace{6pt}

{\bf Acknowledgments.} This research was partially supported by Deutsche Forschungsgemeinschaft through the
SFB 649 Economic Risk.

\end{document}